\documentclass[journal]{IEEEtran}

\ifCLASSINFOpdf
    \usepackage[pdftex]{graphicx}
\else
    \usepackage[dvips]{graphicx}
\fi
\usepackage{url}
\usepackage{paralist}
\usepackage{enumitem} 
\usepackage{color}

\usepackage{amsmath}
\begin{document}
%
\title{A New Take on Protecting Cyclists in Smart Cities}

\author{Adam~Herrmann,
        Mingming~Liu, Francesco Pilla,
        and~Robert~Shorten~
\thanks{A. Herrmann, M. Liu, and R. Shorten are with the School of Electrical, Electronic and Communications Engineering, University College Dublin, Ireland (e-mail: adam.herrmann@ucdconnect.ie; mingming.liu@ucd.ie; robert.shorten@ucd.ie). F. Pilla is with the Dept. of Planning and Environmental Policy, University College Dublin (email: francesco.pilla@ucd.ie).}}

%



\maketitle


\begin{abstract}
	Pollution in urban centres is becoming a major societal problem. While pollution is a concern for all urban dwellers, cyclists are one of the most exposed groups due to their proximity to vehicle tailpipes. Consequently, new solutions are required to help protect citizens, especially cyclists, from the harmful effects of exhaust-gas emissions. In this context, hybrid vehicles (HVs) offer new actuation possibilities that can be exploited in this direction. More specifically, such vehicles when working together as a group, have the ability to dynamically lower the emissions in a given area, thus benefiting citizens, whilst still giving the vehicle owner the flexibility of using an Internal Combustion Engine (ICE). This paper aims to develop an algorithm, that can be deployed in such vehicles, whereby geofences (virtual geographic boundaries) are used to specify areas of low pollution around cyclists. The emissions level inside the geofence is controlled via a coin tossing algorithm to switch the HV motor  into, and out of, electric mode, in a manner that is in some sense optimal. The {\em optimality} criterion is based on how polluting vehicles inside the geofence are, and the expected density of cyclists near each vehicle. The algorithm is triggered once a vehicle detects a cyclist. Implementations are presented, both in simulation, and in a real vehicle, and the system is tested using a Hardware-In-the-Loop (HIL) platform (video provided).
\end{abstract}

\begin{IEEEkeywords} Emissions Regulation; Hybrid Vehicles; Cycling Safety; Smart Cities.
\end{IEEEkeywords}

%
\IEEEpeerreviewmaketitle
\section{Introduction}
Transportation in Ireland and Europe is currently in a state of flux. The last 10 years have seen significant changes in transportation policies, reflecting a move away from emission generating vehicles to more ``greener'' forms of transportation. These include schemes such as: Cycle-to-Work\footnote{http://www.citizensinformation.ie}; low-cost bike rental schemes\footnote{http://www.dublinbikes.ie}; encouragement of hybrid/electric vehicle (HV/EV) ownership through the roll-out of public charging points\footnote{https://www.esb.ie/our-businesses/ecars/charge-point-map}; and tax incentives\footnote{http://www.seai.ie/Grants}. Cycling schemes, in particular, have contributed to the recent increase in Ireland in the popularity of cycling as a means of daily transportation (43\% increase between 2011 and 2016). Such trends are also to be found, not only in Europe, but across the developed world.\newline

With regards to the health \& safety of this increasing cycling population, much focus is naturally placed on the issue of vehicle-cyclist collisions. However, there exists another danger that is often overlooked, and to which cyclists are particularly exposed; namely, the harmful effects of tailpipe emissions from Internal Combustion Engines (ICEs). Many studies exist that document these dangers. According to \cite{nyhan}, in high traffic situations, cyclists are most vulnerable as they are often cycling very close to the exhausts of vehicles who they are sharing the road with. This is made worse by cyclists having a respiratory rate of 2-5 times higher than pedentrians \cite{bigazzi,hoek}. This increased exposure to pollutants contributes to adverse health effects in cyclists populations \cite{cyclingHealth,exposure}. Pollution in our cities has become such an issue that in December 2016 two major European cities, Paris and Athens, have pledged to ban the use of all diesel-powered cars and trucks in their cities by 2025\protect\footnote{https://www.irishtimes.com/news/environment/four-world-capital-cities-to-ban-diesel-vehicles-from-2025-1.2890864}. Car manufacturers, such as Volvo and Volkswagen, have also recognised this problem and have committed to introduce new hybrid and electric vehicles by 2019 and 2020 respectively\protect\footnote{https://spectrum.ieee.org/cars-that-think/transportation/advanced-cars/will-volvo-really-kill-the-gasoline-engine} \protect\footnote{https://www.bloomberg.com/news/articles/2015-09-14/vw-to-reinvent-itself-as-a-maker-of-smartphones-on-wheels}. However, to date, few, if any of these measures are directed at solely cyclists. Our objective in this paper is to address this by proposing simple strategies that can be of immediate benefit to cyclists.\newline

Roughly speaking, traffic emissions in our cities can be addressed in 4 ways: 
\begin{inparaenum}
	\item by building vehicles that don't pollute; 
	\item by informing people of dangerous pollution present in their surroundings so that they can make informed choices; 
	\item by policy interventions such as traffic restrictions (green zones, adaptive speed limits) in certain areas; and  
	\item by using smart devices that adapt to their surroundings to protect humans, such as hybrid actuation in HVs.
\end{inparaenum}
Our approach in this paper is to explore the latter. Due to increased incentives, HVs and EVs have been increasingly adopted by the market over the last a few years. These vehicles have lower emissions (zero emissions in the case of EVs) compared to equivalent pure ICE propelled vehicles. Most HVs are also capable of running on full electric mode without producing emissions, which introduces new ``degrees of freedom in addressing pollution regulation" (see \cite{hybridSwitching}) - especially when groups of vehicles work in tandem. In particular, with these different actuation possibilities, HVs can automatically adapt its drive-train propulsion system in a context-aware manner for the benefit of others - especially cyclists \cite{twitter}. \textcolor{black}{Note that while other strategies can be adopted to reduce emissions along a particular route, for example by managing vehicle speed, or by re-routing vehicles, our approach is to use only the engine actuation mechanisms in hybrid vehicles to lower emissions. Our rationale for doing this is that the method: (i) is completely non-invasive on the driver side; and (ii) does not affect traffic patterns. Thus, it causes the least inconvenience to the driver, especially when compared with informing the driver to lower their speed or to re-route; and avoids disturbances to traffic patters due to re-routing of traffic or reduced flow rates. Note also that strategies to reduce emissions via speed adjustment, or by rerouting, is reported in our previous work; see \cite{Cri17} and the references therein. } \newline

In this paper, the actuation possibilities of HVs are used to control the emissions in a virtual geographic boundary (geofence), whose coordinates are determined by the location of the vehicle that has detected a cyclist, using Radio Frequency IDentification (RFID) technology. The system presented is one aspect of a larger body of our current work which aims to deal with multiple aspects of cycling safety issues encompassing:
\begin{inparaenum}
	\item cyclist detection for a vehicle driver alerting system using RFID technology;
	\item enhancing these alerts with coarse localisation of the cyclist's position relative to the vehicle;
	\item situational-aware pedal assisted e-bike to regulate the breathing rate of the cyclist in high pollution areas \cite{shaun}.
\end{inparaenum}
The details of these systems will not be discussed as their implementations are outside the scope of this paper; however, several aspects are used to implement the solution proposed. Further, we note that this paper is dealing solely with \textit{pollution} and the harmful effects caused to humans by high concentrations of pollutants in their environment. This topic is \textit{not} related to greenhouse gas effects. \newline

The remainder of this paper is organised as follows. Previous works are reviewed in Section II. Description of the system model and the suggested operation modes are given in Section III. Simulation set-up in Simulation of Urban MObility (SUMO) and the details of  Hardware-In-the-Loop (HIL) implementation are presented in Section IV. Simulation results are discussed in Section V. Finally, future work and our conclusion are presented in Section VI and VII, respectively. 


\section{Previous Work}
Our RFID cyclist detection system builds on the work of \cite{colum} which introduced the design and successful implementation of a system to detect and alert the driver of a nearby cyclist using passive RFID tags on the cyclist and an RFID antenna located on the vehicle. This work enables the detection of the cyclist to trigger the pollution mitigation system. Our work is also related to the wide body of literature on engine management systems for hybrid electric vehicles. In \cite{instantaneous} and \cite{shev}, engine management systems (EMS) were developed to optimally control the energy usage of a HV to minimise fuel consumption. In \cite{overview}, a detailed overview of both parallel and series HVs and their power management strategies were presented. In \cite{indirectly}, an EMS was developed to minimise air pollution caused by the HV both directly and indirectly (whilst driving as well as vehicle/grid interaction). Furthermore, in \cite{hybridSwitching} \cite{Mahsa} and \cite{arieh} an emission trading framework was proposed whereby the problem of sharing an emission budget between HVs is formulated as a utility maximisation problem. The authors in \cite{arieh} used this framework to further implement a pollution mitigation system in a real HV. This system works by using a GPS-enabled Android smartphone to automatically switch the HV into EV Mode when inside fixed geofences around residential or school areas. In \cite{twitter} an EMS was formulated which aims to lower the emission output of the vehicle in areas of high pedestrian traffic taking account of the uncertainty of routes the driver may travel to. Our paper will extend the work of \cite{arieh} by using the RFID system to create the geofence dynamically and have it move with the cyclist rather than remain static. Our work will also use the ideas from \cite{twitter} by incorporating popular cyclist routing data to lower emissions on roads in the geofence which cyclists use often, thereby predicting the route of the cyclist which has been detected in the geofence. 

\section{System Model and Operations}
\label{probsol}

Our objective is to design a system to lower the  pollution in the area around cyclists to some level; in other words to attempt to create a bubble of clean air around a single or group of cyclists. In particular, we shall discuss two modes of system operation; namely a single-vehicle mode, and a multi-vehicle mode, depending on the radius of the area around the cyclist that is to be protected.

\subsection{Basic Set-up}
\begin{figure}
    \centering
    \includegraphics[width=3.3in]{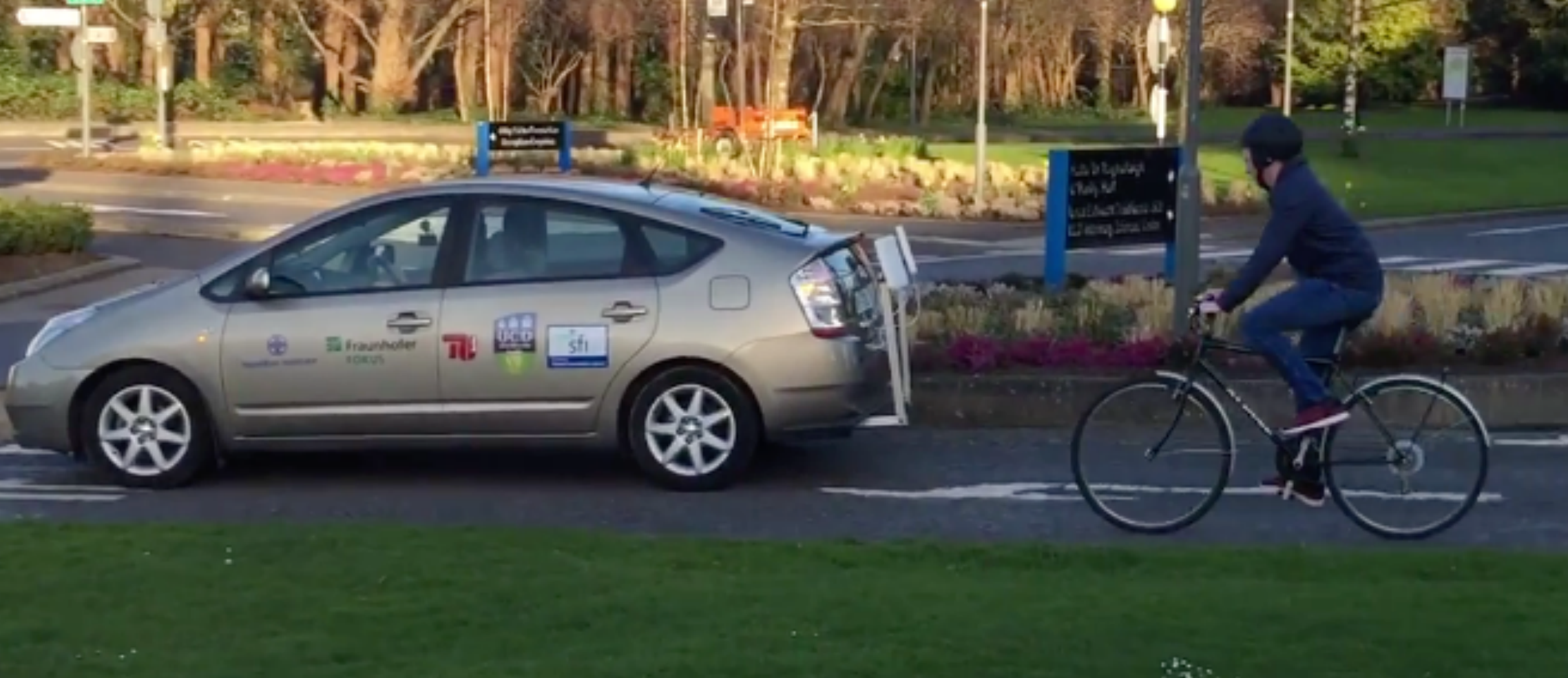}
    \caption{Toyota Prius with the RFID antennas mounted on the rear. A passive RFID tag is mounted on the cyclist's helmet.}
    \label{prius}
\end{figure}

The basic idea is very simple. This system uses the actuation abilities of HVs to dynamically switch the vehicle from a ``polluting'' mode to an emission free mode (i.e. full electric mode) when a cyclist is detected. This fully automated system requires no driver interaction as the transition between polluting and electric modes are designed to be seamless. Specifically, a cyclist is detected by the RFID cyclist detection system developed by \cite{colum}. The detection system uses a RFID reader mounted on the rear of a real HV (Toyota Prius - Fig. \ref{prius}). The reader can detect the presence of an RFID tag mounted on the cyclist's helmet. A Java Server located in the rear of the vehicle processes the information coming from the RFID readers and sends the relevant information to an Android application to inform the driver whether there is a cyclist nearby and what status the system is in. Given this info, the Prius switches between electric and polluting modes automatically (without driver's interaction) through a Bluetooth Actuator developed in \cite{arieh}. We note that the cyclist detection system could be replaced with image recognition using cameras or radar mounted on the vehicle, as discussed in \cite{camera}. Our preference for RFID is due to its passive nature (no battery required) and low cost on the cyclist side, the reduced privacy concerns when compared with cameras, and the low computational cost of the technology on the vehicle side. 

\subsection{Single-Vehicle Operation}

The single-vehicle operation mode is designed to simply switch the HV to full electric mode when  a cyclist is detected by the vehicle, and for a period thereafter. This solution is designed to provide protection to a cyclist who is cycling directly behind (or in proximity to) the vehicle by eliminating tailpipe emissions\footnote{A demo video for this test can be seen at \url{https://youtu.be/8Qnoc5GD7kI}}. While this system always reduces the direct impact on the cyclist approaching to the vehicle, and protects the cyclist from emission spikes from tailpipes, it does not take into account many other factors that can affect cyclists; for example, other polluting vehicles nearby; other sources of pollution; and the effects of route choices that a cyclist may make. For these reasons, we consider a generalisation involving multiple vehicles in the following section.\newline

\subsection{Multi-Vehicle Operation (Geofences)}

In this section, we introduce a multi-vehicle operation mode that can be deployed by vehicles inside the geofence to protect the detected cyclist. The suggested system is operated as follows. When a cyclist is detected, the vehicle which detects the cyclist will broadcast a signal to nearby vehicles (as specified by the geofence) to switch on the system. Vehicles receiving the signal (vehicles inside the geofence) can proactively switch into electric engine mode in a way that the aggregate emission from the group inside the geofence are controlled at a given level. Note that the total damage done to the cyclist (by the vehicles) is a function of his/her breathing rate, the average aggregate emission level, and the exposure time. Switching in and out of electric mode is done at regular intervals. The switching decision of each vehicle in the geofence is calculated based on the polluting level of the vehicle, and the likelihood the cyclist will encounter this vehicle as he/she travels. In the current embodiment of the system, this likelihood is estimated using historical data. This is available from local authorities,  fitness tracking apps, and other sources. We use data from a fitness application for the proof of concept in this paper. In what follows we make the following assumptions. \newline
\begin{enumerate}
	\item First, we assume that a number of vehicles in the geofence are either HVs or EVs. Note that this system could also work for traditional ICE vehicles by altering the engine characteristics. Vehicle manufactures currently offer different ICE operating modes to drivers such as Eco and Sport.\newline
	\item When a vehicle is in electric mode, it is emitting no pollution and thus causing no harm. When the vehicle is in Polluting Mode, it is emitting pollution characterised by its engine class.\newline
	\item The HVs are parellel hybrid vehicles (eg. Toyota Prius), which allow different modes of operation (full electric or hybrid). \newline
\end{enumerate}

\noindent\textbf{Remark: } Our prototype system is currently not designed to take into account the distance from individual vehicles to cyclists. This stems from the fact that it is difficult to estimate distance from RFID RSSI (Return Signal Strength Indicator) data due to the short range and environmental influence on the data values. However, since this system deals with lowering the aggregate emissions in the geofence, this is not an issue. Our system could of course be extended to account for distance by instrumenting the vehicles with the appropriate technologies. It is, of course, possible to force the detecting vehicle to switch to ECO mode as this is most likely the vehicle closest to the cyclist.  

\begin{figure}
    \centering
    \includegraphics[width=3.3in]{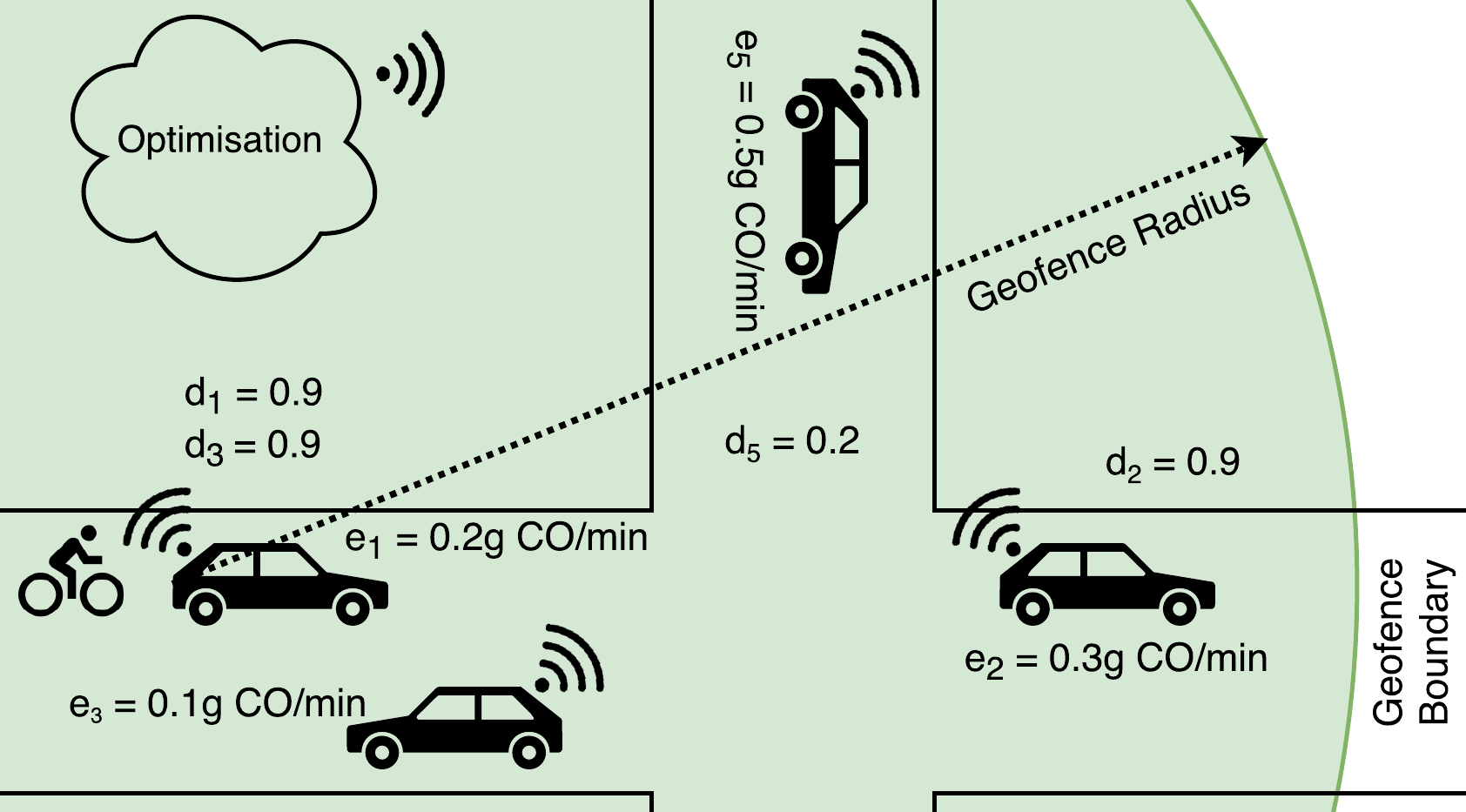}
    \caption{Example map showing cyclist being detected and geofence being created. The values of $d$, the estimated cyclist density value, applies to all cars on the road edge. Each vehicle has its own CO emission value, $e$ (grams of CO/min) calculated using the average speed model and is based on the euro class rating.}
    \label{geofence_ill}
\end{figure}

\subsection{Optimisation}

\textcolor{black}{For a fixed  level of group emissions the operation of all vehicles in the geofence is characterised by the following simple optimisation and is illustrated in Fig. \ref{geofence_ill},}
 
\begin{equation}
\begin{aligned}
& \underset{x_i}{\text{max}}
& & \sum_{i=1}^N \frac{x_i}{d_i}\\
& \text{s.t.}
& & \sum_{i=1}^N x_ie_i \leq {E}(\Delta)\\
& & & 0 \leq x_i \leq 1
\end{aligned}
\label{eq1}
\end{equation}
\textcolor{black}{where $x_i$ is the decision variable (for example, a probability that controls the likelyhood that the vehicle $i$ is in polluting mode at a given instant); $d_i$ is the estimated cyclist density value for the road segment where the $i$'th vehicle is on; $e_i$ is the estimated CO emission (gCO/min) of the $i$'th vehicle; $E(\Delta)$ is an emission limit for the cars in the geofence (gCO/min) at the decision instant; $\Delta$ is a measure of the difference between the measured pollution concentration and a safe level (see remark below); and $N$ is the total number of cars in the geofence. The optimisation step of the algorithm is repeated at an interval of $\tau$ seconds where $E(\Delta)$ may vary from one decision to the next (again see below). Note that in our embodiment of the system, $x_i \in [0,1]$ can be chosen to be either the fraction of time between decisions that a vehicle is in the polluting mode, or can be chosen to be a probability that vehicle $i$ remains in polluting mode for an interval of time. In our results section that follows, this latter interpretation is implemented both for ease of execution and to be consistent with our prior work \cite{Cri17}. 
It can be observed that the vehicles that are travelling on roads having lower $d_i$ and $e_i$ are favoured (i.e. more likely to stay in the polluting mode).}\newline 

\noindent\textcolor{black}{{\bf Remark (Mode Switching):}  The interpretation of $x_i$ as probabilities is sometimes preferable as it may achieve a multiplexing effect and avoid pollution spikes inside the geofence. Here the idea is that every $t$ seconds, each HV changes mode and switches into polluting mode with probability $x_i$. Similar probabilistic switching strategies are discussed in \cite{hybridSwitching, Cri17, Sponge}.}\newline

\noindent{\textcolor{black}{{\bf Remark (Background Emissions\protect\footnote{The authors gratefully acknowledge the comments of the anonymous reviewer for pointing out ambiguities in the original manuscript.})}: Before proceeding, it is worth noting that not all emissions are as a result of transportation effects. Background emissions are always present and contribute to the measurable aggregate emissions. The objective of our strategy is to regulate the aggregate emissions inside a geo-fence by managing the switching of a network of hybrid vehicles into their polluting mode in an elastic manner. One approach to do this is to embed our optimisation strategy inside of a feedback loop as depicted in Fig.~\ref{feedback}. 
\begin{figure}
    \centering
    \includegraphics[width=3.3in]{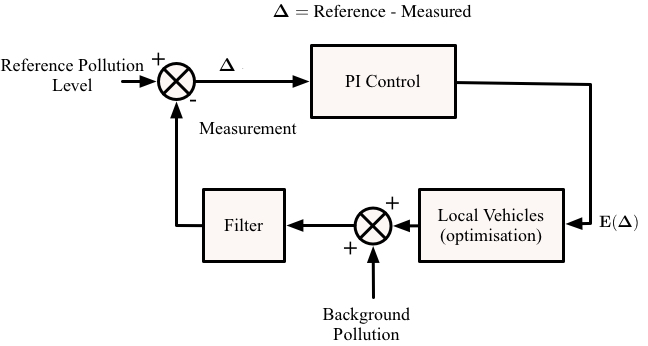}
    \caption{Feedback strategy, background emissions and disturbance rejection}
    \label{feedback}
\end{figure}
Here, the feedback loop regulates the {\em aggregate emissions} by determining the rate limit $E(\Delta)$. This is done by comparing the measured aggregate with a fixed limit as specified, for example, by EC directive $2008/50/EC$\protect\footnote{http://ec.europa.eu/environment/air/quality/standards.htm} \protect\footnote{http://eur-lex.europa.eu/legal-content/EN/TXT/?uri=celex{\%}3A32008L0050}, or a world health organisation (WHO) equivalent\protect\footnote{http://www.euro.who.int/\_\_data/assets/pdf\_file/0005/74732/E71922.pdf}, as part of a feedback strategy. The role of the optimisation is then to manage the elastic component of the overall strategy; namely by adjusting the overall local emissions of the group of hybrid vehicles. This approach is akin to treating the non-measurable part of the emission as a disturbance, and rejecting the disturbance by varying $E(\Delta)$. A perhaps simpler alternative is to determine $E(\Delta)$ by referencing the allowable aggregate emissions to the {\em urban background} as defined as per EC Directive 2008/50/EC. The urban background pollution (CO, NOx, PM, etc) could be retrieved in real-time from official monitoring stations around the city or using low cost sensing. In both cases, adjusting $E(\Delta)$ seems reasonable due to the emergence of low cost pollution sensing\protect\footnote{http://iscapeproject.eu/} in our cities \cite{francesco}. It is also worth noting that, typically, the background emissions will vary slowly, both spatially and temporally so that $E(\Delta)$ will  not vary too much from one decision instant to the next.  Notwithstanding this comment, more dynamically varying background levels of pollution can be incorporated into our framework by decreasing the time between optimization decisions.  Note also that in the case that $E(\Delta) \leq 0$  all cars switch to non-polluting mode, and help the cyclist even when background pollution levels are high. Finally, note that the idea of embedding a group of hybrid vehicles inside a feedback loop has been implemented in the context of our previous work \cite{hybridSwitching, Cri17} and we refer the interested reader to these papers for more details.}\newline

\noindent{\textcolor{black} {{\bf Remark (Prioritisation of Cyclists)} : In some circumstances, the aforementioned strategy may have the effect of prioritising cyclists over pedestrians and other citizens. The justification for this, in the present context,  is that cyclists are typically more vulnerable than other road users due to their increased respiratory rate, and their proximity to vehicle tailpipes (with some studies suggesting that they experience more than twice the level of exposure than pedestrians). However, the reader should also note that additional constraints can also be introduced into the optimisation to account for other road users and city dwellers; see \cite{hybridSwitching, twitter, Cri17}.}}\newline

\subsection{Cyclist Density Data}

Cyclist density data required for the optimisation  \eqref{eq1} was obtained from historical data provided by  the cycling route mapping service Strava\footnote{http://labs.strava.com/heatmap/}. The data used in this paper is shown in Fig. \ref{fig_strava} where we assigned each road with a weight based on the heatmap data from Strava.
\begin{figure}
    \centering
    \includegraphics[width=3.3in]{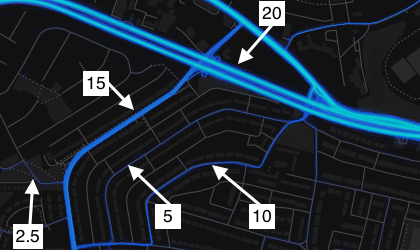}
    \caption{Heat-map showing previous routes of cyclists in the Palmerstown area. Brighter lines indicate more popular routes. The relative weights assigned for use with the optimisation are also shown overlaid. If a road has no data it is assigned a value of 1.0 by default.}
    \label{fig_strava}
\end{figure}
Note that more accurate and real time data can be obtained from sources such as local authorities, mobile phone companies, and social networks.


\subsection{Emission Data}
\textcolor{black}{For our simulations, and our hardware-in-the-loop emulation (see below),} the pollution data for each car in the geofence is estimated using the average speed model proposed in \cite{hybridSwitching}. The emission function (g/km) for the vehicle in question, where $t$, $v$, and $p$ are the type of vehicle, average speed and pollutant respectively, is repeated in \eqref{em}:
\begin{equation}
f(t,p)=\frac{k}{v}(a+bv+cv^3+dv^3+ev^4+fv^5+gv^6)
\label{em}
\end{equation}
where the parameters $k, a, b, ..., g$ are used to specifiy different level of emissions by different classes of vehicles. The functions are converted to pollutants per minute and the aggregate is calculated by summing all vehicles in polluting mode. \newline

\noindent\textcolor{black}{{\bf Remarks on the {\em average speed} model}: We make the following comments concerning the use of the {\em average speed} model.}} \newline

\begin{itemize}
\item \textcolor{black}{The average speed model is used in this paper to illustrate the operation of the system, and to facilitate our SUMO based simulations and hardware-in-the-loop testing. 
It is {\bf NOT} intended to be a realistic representation of vehicle emissions.}\newline

\item \textcolor{black}{In practice, measurement data from each vehicle type, and a {\em Portable Emissions Management System} (PEMS),  would be used to construct a 
more accurate emissions model\protect\footnote{https://ec.europa.eu/jrc/en/vela/portable-emissions-measurement-systems} \protect\footnote{http://www.horiba.com/automotive-test-systems/products/emission-measurement-systems/on-board-systems/}\cite{Gie16}. Thus, in practice, a PEMS, could be used to calibrate the algorithms to model emissions for different driving cycles over an inventory of vehicles.
Alternatively, emissions could be measured in an active manner using PEMS like devices.}
\newline

\item \textcolor{black}{We note again, that control of emissions, via adjustment of vehicle speed, is {\bf NOT} the subject of this paper. We are interested in non-invasive strategies (from the viewpoint 
of the driver, and traffic disturbance). Consequently, control of switching between polluting and non-polluting mode in the HV, is our preferred actuation method. Finally, as we have already noted, 
pollution mitigation measures, such as speed control, or traffic rerouting, have been investigated in a series of our companion papers \cite{Cri17}.} 
\end{itemize}

\subsection{Implementation Aspects for Multiple Vehicle Simulation and Emulation}

Finally, we comment on implementation aspects of the multi-vehicle operation mode.\newline 

\begin{enumerate}
    \item Vehicles begin in polluting mode. When a cyclist is detected by the RFID system, an alert is sent to the driver and the central server, which creates the geofence with a unique ID corresponding to the ID of the RFID tag on the cyclist's helmet. The centre of the geofence will be assigned to the location of the vehicle that detected the cyclist.\newline  
    \item All vehicles inside the geofence are asked for their emission class (how polluting they are) and what road they are on. The cyclist density value for the road the vehicle is travelling on is obtained from cycling statistics data for the area, and updated periodically.\newline
    \item The central server solves the linear programming problem \eqref{eq1} and calculates the probability of staying in polluting mode (not switching into electric mode) for each vehicle, based on how polluting the vehicle is and how likely it is that a cyclist will travel on the road the vehicle is on.\newline
    \item These probabilities are communicated to each vehicle who act on this instruction by performing a weighted coin toss on whether or not the vehicle will switch into electric mode.\newline
    \item Steps 2 to 4 are repeated every $\tau$ seconds until the cyclist is no longer detected by a vehicle in the geofence. If a cyclist has not been detected for 20 seconds, the geofence is removed and the vehicles inside the removed geofence switch back to polluting mode.\newline
\end{enumerate}

\section{SUMO and HIL Simulation}

In this section we evaluate our proposed method using the road-network and traffic simulator, SUMO, and in a real vehicle using the Hardware-In-the-Loop (HIL) platform. Note SUMO \cite{wynita} is an open source, microscopic road traffic simulation package from the Institute of Transportation Systems at the German Aerospace Centre (DLR).

\subsection{SUMO Set-up}

\begin{figure}
	\centering
	\includegraphics[width=3.3in]{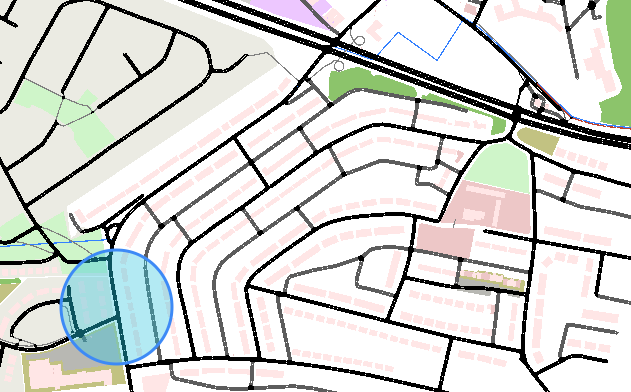}
	\caption{Map in SUMO obtained from Open Street Map of Palmerstown, Dublin 20.}
	\label{sumo_palmerstown}
\end{figure}

A SUMO simulation was conducted using a map of the Palmerstown, Dublin 20, in Ireland obtained from Open Street Map\footnote{https://www.openstreetmap.org}. The heatmap data from Strava shown in Fig. \ref{fig_strava} applies to the map shown in Fig. \ref{sumo_palmerstown}. The Palmerstown area was used for the simulation as the heatmap from Strava has a diverse range of cyclist density values in a small area, which makes it ideal for testing the geofence algorithm. In the simulation, vehicles are added at every time step with each vehicle being assinged a EURO emission class between 1 and 4\footnote{The EURO emission classes define limits for exhaust-gas emissions. A higher emission class corresponds to stricter emission limits.}. The vehicles will be removed from the road network when they arrive at their pre-defined destinations. In particular, we assume that a cyclist will be detected by vehicle 24 at a specific time step, and a geofence centred at the vehicle 24 will be created correspondingly \footnote{In the case of a real vehicle, a geofence will be created automatically when the RFID system detects a cyclist around.}. The optimisation algorithm is implemented at each time step for the vehicles inside the geofence and the vehicles act on this information as soon as the results of the optimisation are received from the central server. Furthermore, the emission limit for all vehicles inside the geofence is set to 1.0gCO/min and the radius of the geofence is set to 100m. These limit values are chosen arbitrarily to demonstrate the  proof of concept. 

\subsection{HIL Set-up}

In order to get a real-world feel for the system developed, our algorithm was implemented and tested in a real car using the HIL platform. This platform allows a real vehicle to control an avatar vehicle inside the SUMO simulation. This creates the effect of having a fleet of connected vehicles while only a single vehicle is needed in reality. A full description of the HIL platform can be seen in \cite{wynita}; a brief overview follows.\newline

\begin{figure}
	\centering
	\includegraphics[width=0.9\linewidth]{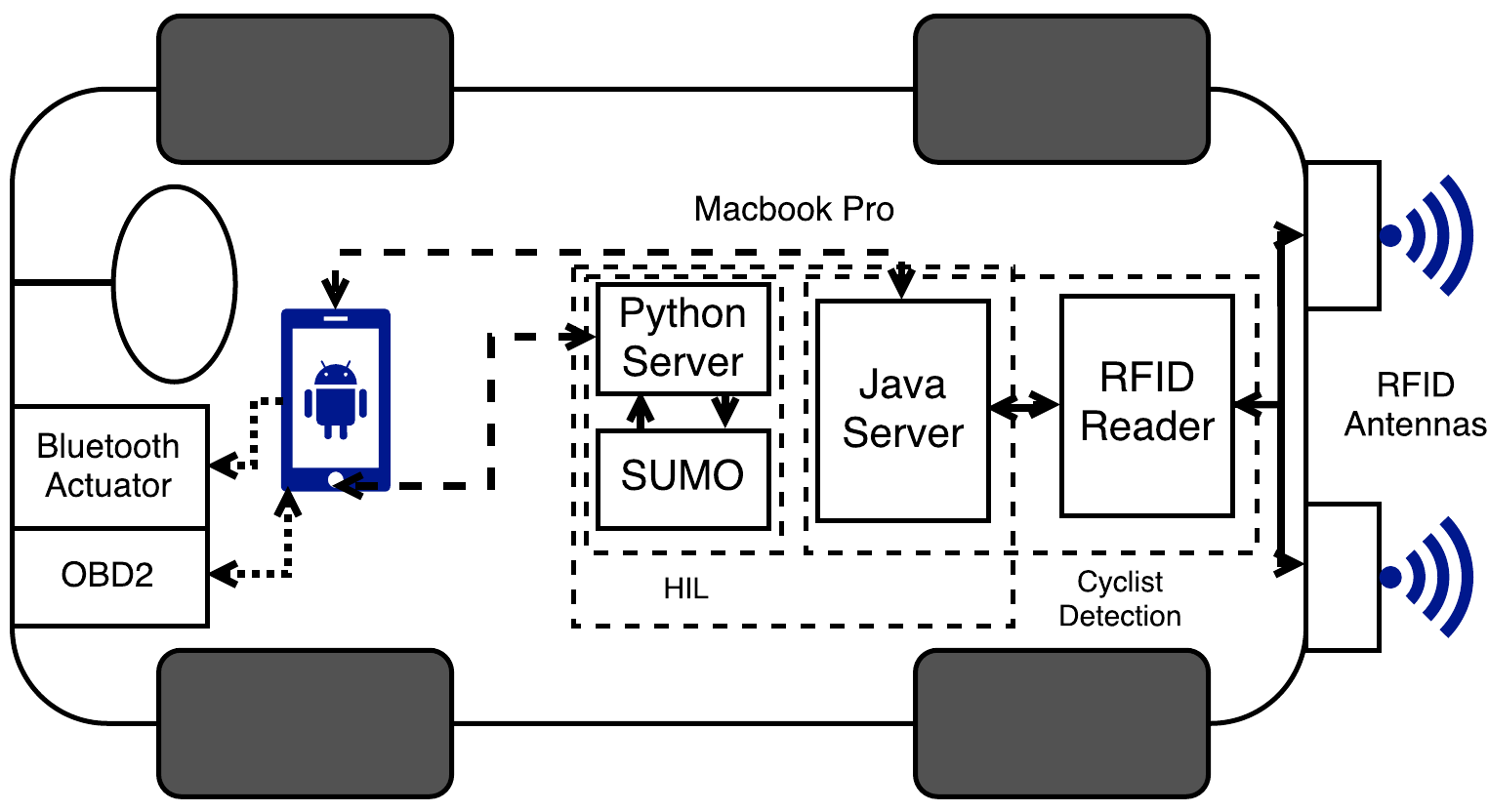}
	\caption{In-car implementation of the system. Solid lines indicate wired ethernet connection, dashed lines indicate WiFi connection and dotted lines indicate Bluetooth connection.}
	\label{hil}
\end{figure}

The system architecture used for this work can be seen in Fig. \ref{hil}. For logistical reasons, the UCD campus was used for testing rather than Palmerstown. However, the data from Strava is not as diverse for the UCD campus as it is a very busy cycling location as many students cycle to and from campus. \\

\noindent \textbf{Vehicle:} The vehicle and cyclist detection system used for the HIL test is similar in configuration to the one described in Section \ref{probsol}. The vehicle controls the avatar vehicle by sending the current speed of the vehicle to the SUMO Simulation. The speed is obtained from the OBD2 port of the vehicle via a Bluetooth OBD2 Module connected to the Android application via the Torque Pro Libraries\footnote{https://torque-bhp.com}. The Android application handles communication with the Bluetooth OBD2 module and relays the speed information to the Python Server (running the SUMO Simulation) via WiFi.\\

\noindent \textbf{On-Board Computer:} A Macbook Pro (15" Mid 2009) runs the Java and Python servers for the RFID System and SUMO simulation respectively. These two servers act as two individual components of the system and do not communicate on the computer (to simulate the Python server running in the cloud, while the Java server runs on the vehicle). The two servers communicate through the Android application.

\subsection{HIL Implementations}

When a cyclist is detected by the RFID reader, the Python server is instructed to create/update a geofence at the location of the real vehicle in the simulation. The results of the optimisation are then communicated with the vehicles in the geofence, if the real vehicle is one of those vehicles, the result as to whether to switch on electric mode is sent to the Android application which instructs the Bluetooth actuator to turn on electric mode. Since we are dealing with a real vehicle, a couple of seconds are required to execute the ``change to electric mode'' command. As such the frequency that the optimisation runs inside the geofence was changed from every second (as in the SUMO simulation), to every five seconds to allow for time for the vehicle act on the commands.\\

\section{Results and Discussion}

In this section we  discuss our simulation results. In what follows we assume that the emission rate limit $E(\Delta)$ is constant; this is akin to assuming that the background emissions are constant.  A video of the SUMO simulation and the HIL test, illustrating the driver warning system and the automatic operation of the pollution mitigation system can be seen at \url{https://youtu.be/wWU7kREc_2A}.\newline

\begin{figure}
   \centering
   \includegraphics[width=0.9\linewidth]{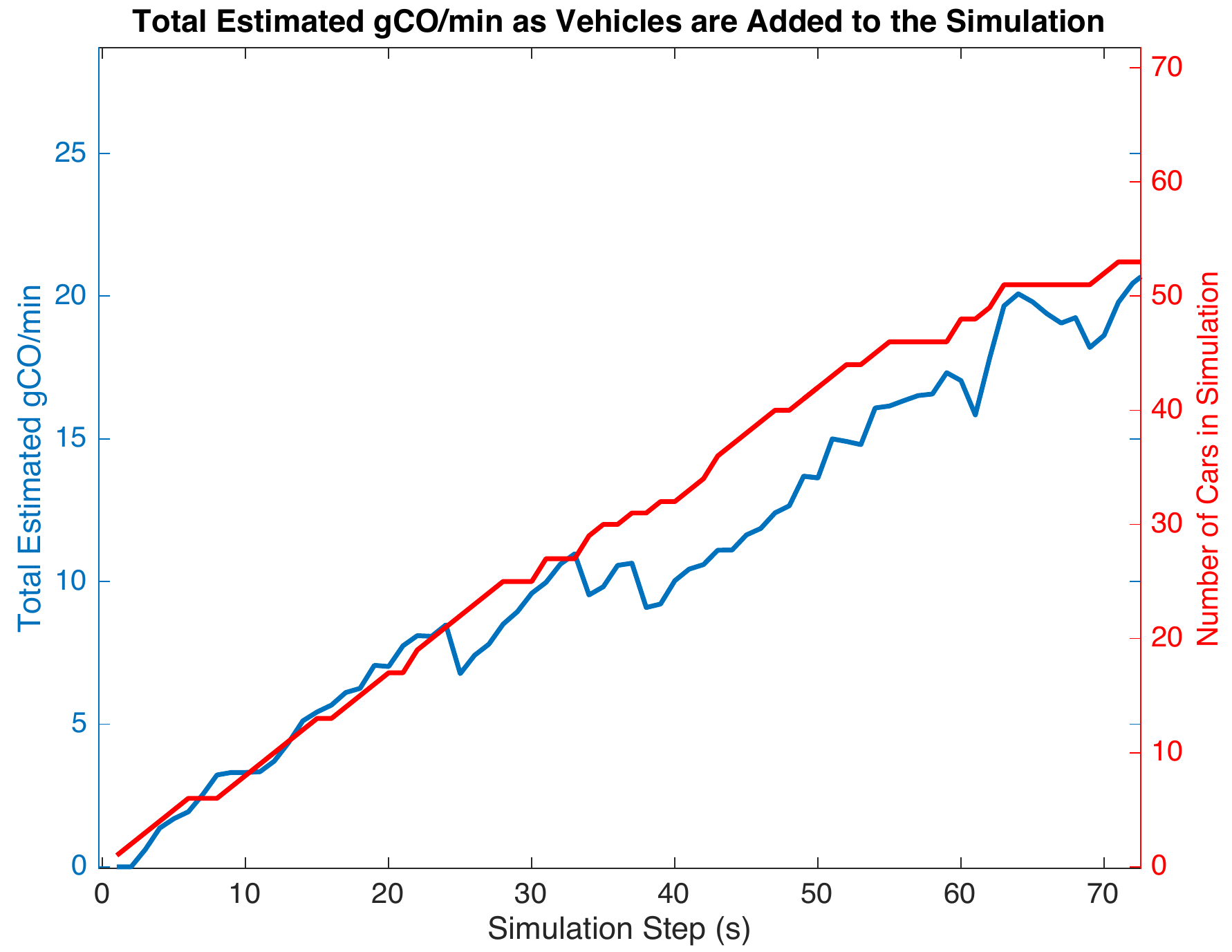}
   \caption{As more vehicles get added to the SUMO Simulation the total emissions increase.}
   \label{geofence_total}
\end{figure}

\begin{figure}
    \centering
    \includegraphics[width=0.9\linewidth]{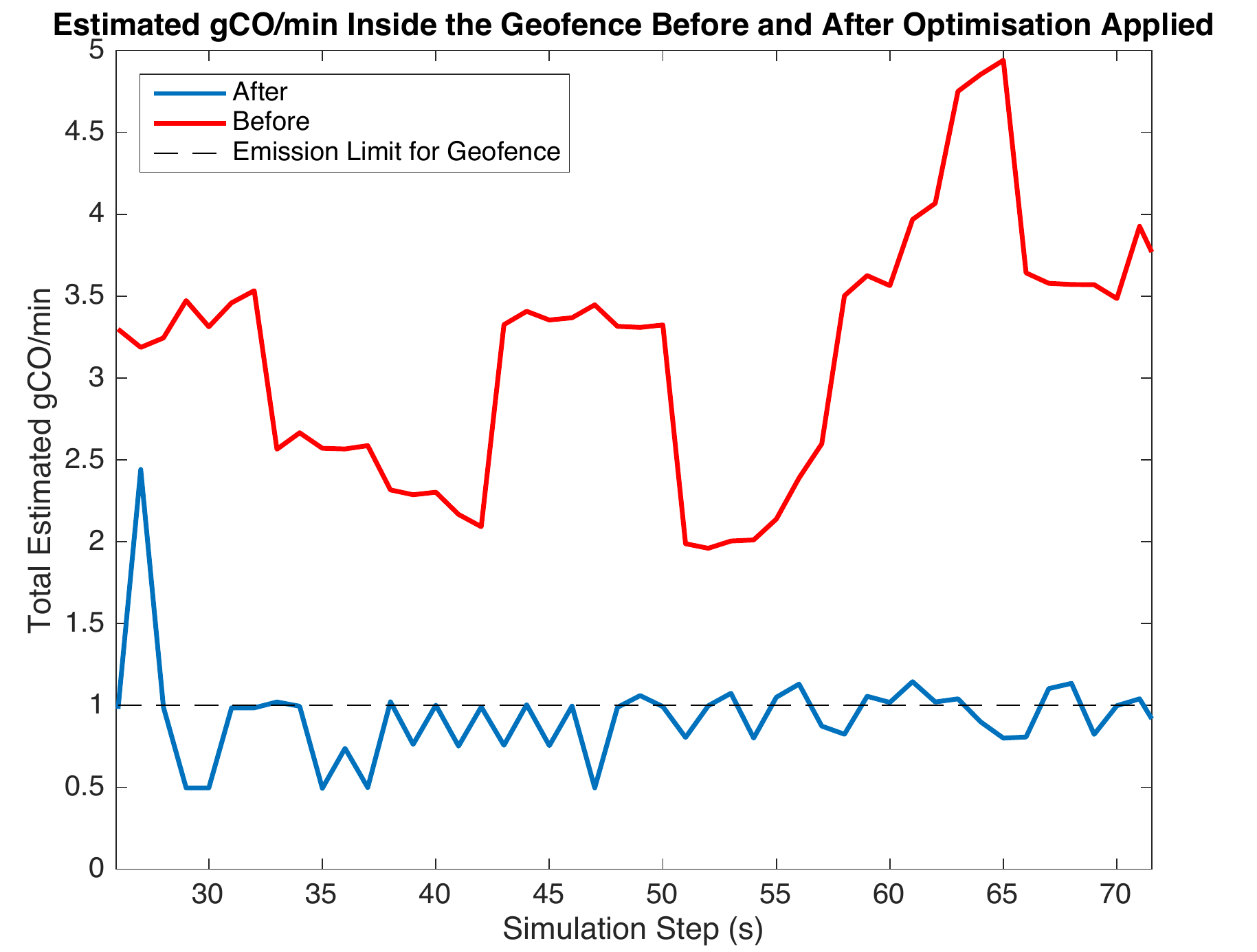}
    \caption{Total emissions inside the geofence before and after the algorithm is applied for each time step.}
    \label{geofence_beforeafter}
\end{figure}

\begin{figure}
    \centering
    \includegraphics[width=3.3in]{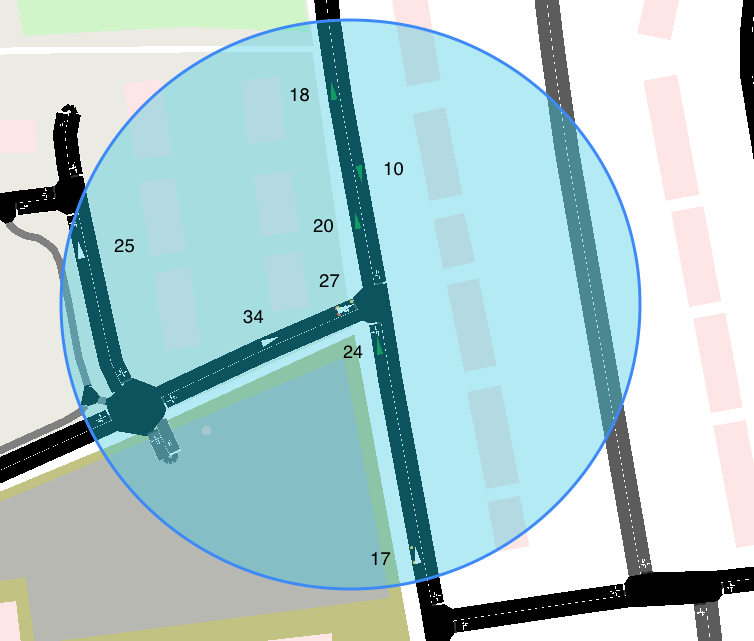}
    \caption{Geofence operating in SUMO. Green vehicles indicate vehicles currently in electric Mode, where Grey vehicles are in Polluting Mode. The cyclist was detected by the vehicle in the centre of the geofence.}
    \label{geofence_overlay}
\end{figure}

Fig. \ref{geofence_total} shows the total emissions for the map as more vehicles are being added. Fig. \ref{geofence_beforeafter} shows the results that this pollution mitigation technique has achieved. The maximum emission level of $E(\Delta) = 1.0gCO/min$ inside the geofence has been achieved with only slight variations caused by the uncertainty of the nature of the probabilistic assignment of whether the vehicles will stay in polluting mode.
Fig. \ref{geofence_overlay} shows the geofence in operation at a specific time step in SUMO with vehicles 10, 18, 20, and 24 operating in electric mode (green) and vehicles 17, 25, 27, and 34 operating in polluting mode (grey). As expected, more vehicles are in electric mode on the road that has a high cyclist density value. \newline

\begin{figure}
    \centering
    \includegraphics[width=0.9\linewidth]{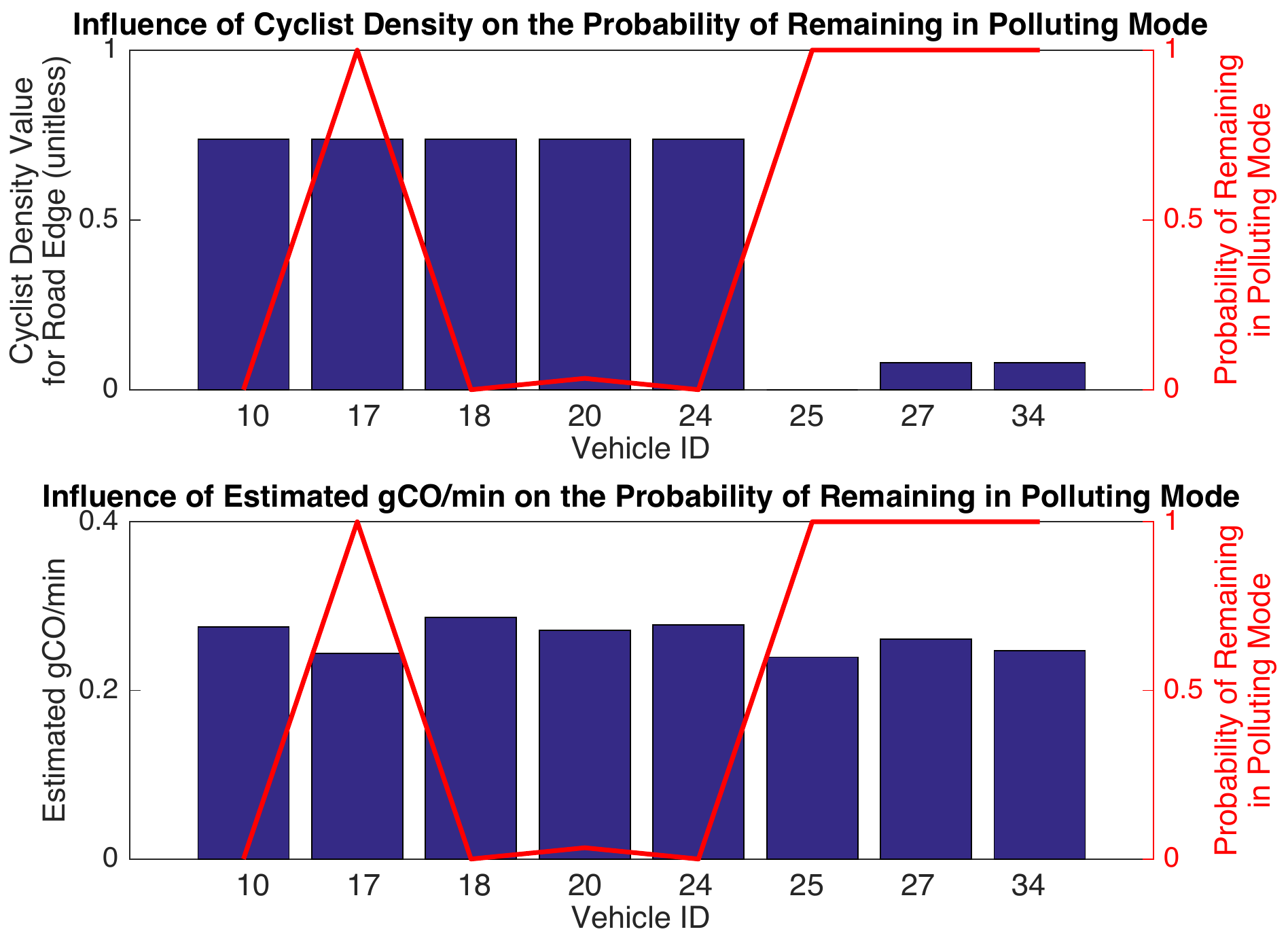}
    \caption{A detailed look at the effect of the emission class of the vehicle and the estimated density of a cyclist taking the road the vehicle is on has on the probability of remaining in polluting mode.}
    \label{geofence_more_cars}
\end{figure}

\begin{figure}
    \centering
    \includegraphics[width=0.9\linewidth]{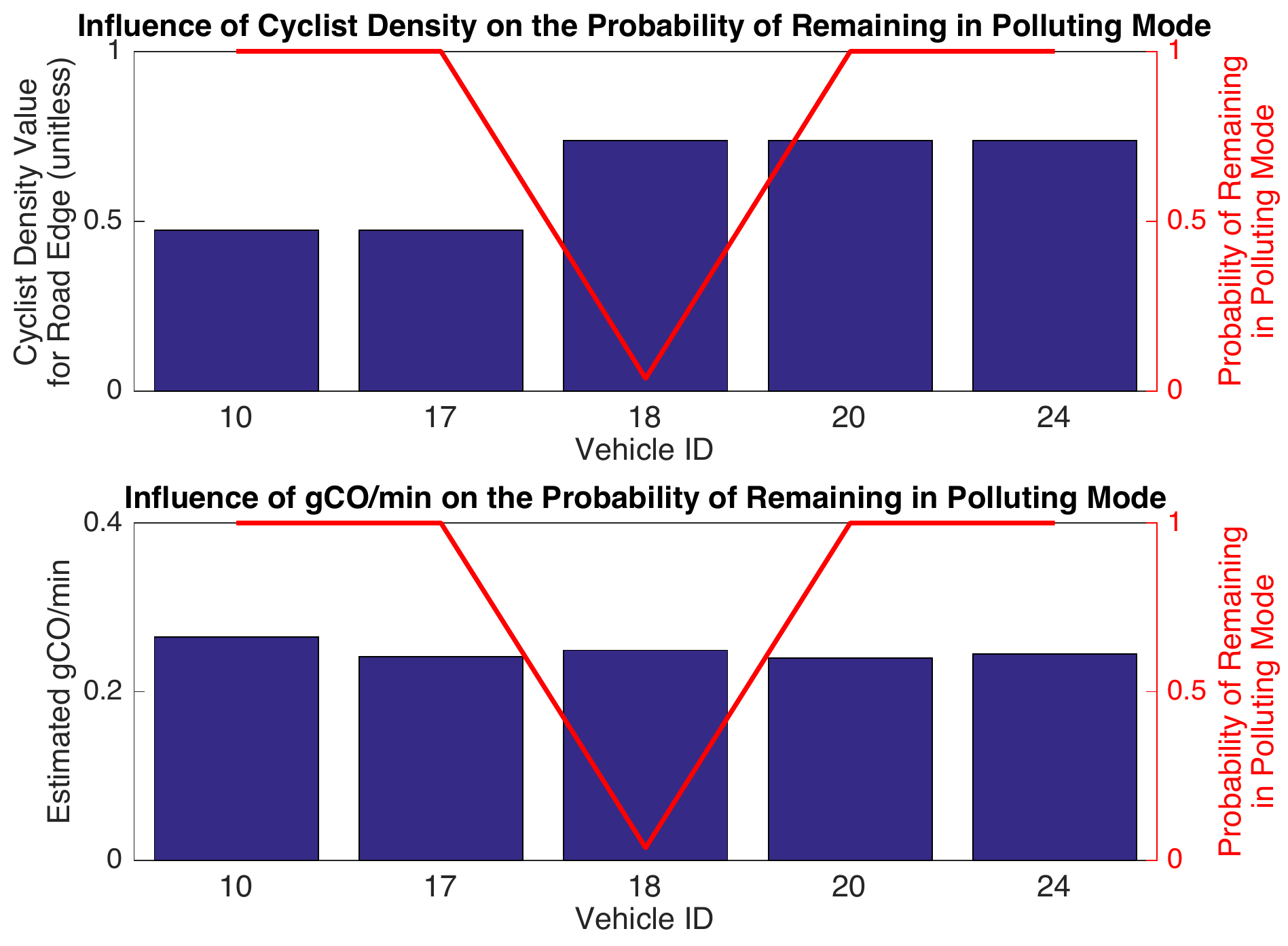}
    \caption{The effect of fewer cars in the geofence.}
    \label{geofence_less_cars}
\end{figure}

Furthermore, Fig. \ref{geofence_more_cars} shows a detailed look at the data the optimisation used and the results assigned to each vehicle. Vehicles 10, 17, 18, 20, and 24 are all on the main road which have a high estimated cyclist density value relative to vehicles 25, 27, and 34 who are on neighbouring roads that are much less likely for a cyclist to travel on, causing them to have a high probability of remaining in polluting mode. Here the emission constraint has been satisfied by the vehicles on roads of high estimated cyclist density switching to electric mode; however, vehicle 17 has been selected to remain in polluting mode despite being on a road of high estimated cyclist density. This is due to: 
\begin{inparaenum}
    \item The emission levels being safe enough to allow vehicle 17 to remain in polluting mode (the utility of the EVs limited battery capacity is being maximised); and
    \item Vehicle 17 not being as polluting compared to the other vehicles on the road and in the geofence, so if there is spare `utility', vehicle 17 is preferred.
\end{inparaenum}
Therefore, this vehicle stays operating in polluting mode as seen in Fig. \ref{geofence_overlay}.
Finally, Fig. \ref{geofence_less_cars} shows how fewer cars in the geofence effect the results. Since there are fewer cars in the geofence, a lower quantity of overall emissions are being generated, so the optimisation can allow more vehicles to travel in polluting mode. Here vehicles 10, 17, 20, and 24 are allowed to stay in polluting mode as there is spare `utility'. Vehicle 18 is assigned a high probability of switching to electric mode as out of the three vehicles on the road with high estimated cyclist density (vehicles 18, 20, and 24), it has the highest output of gCO/min.

\section{Future Work}
Future work could involve further development of the geofence algorithm to address issues not investigated in this paper such as: 
\begin{inparaenum}
\item allowing the operation and simulation of multiple overlapping geofences within an area and the impact on the emission sharing; and
\item incorporating the vehicle route prediction ideas from \cite{twitter} to better determine where pollution hot spots are going to occur and stop them \textit{before} they happen, rather than instantaneously.
\end{inparaenum}

\section{Conclusion}
In this paper a novel pollution mitigation technique has been introduced with the focus to protecting cyclists from unsafe quantities of tail-pipe emissions. The technique has been tested both in simulation and in a real vehicle and shows its ability to optimally regulate emissions inside of a dynamically generated geofence around a cyclist. 

\section * {Acknowledgements } The authors gratefully acknowledge the comments of Emanuele Crisostomi and Giovanni Russo. This work was in part supported by Science Foundation Ireland grant 11/PI/1177 and by iSCAPE (Improving Smart Control of Air Pollution in Europe) project, which is funded by the European Community's H2020 Programme (H2020-SC5-04-2015) under the Grant Agreement No. 689954

\bibliographystyle{ieeetran}
\bibliography{References}
\newpage
\begin{IEEEbiography}[{\includegraphics[width=1in,height=1.25in,clip,keepaspectratio]{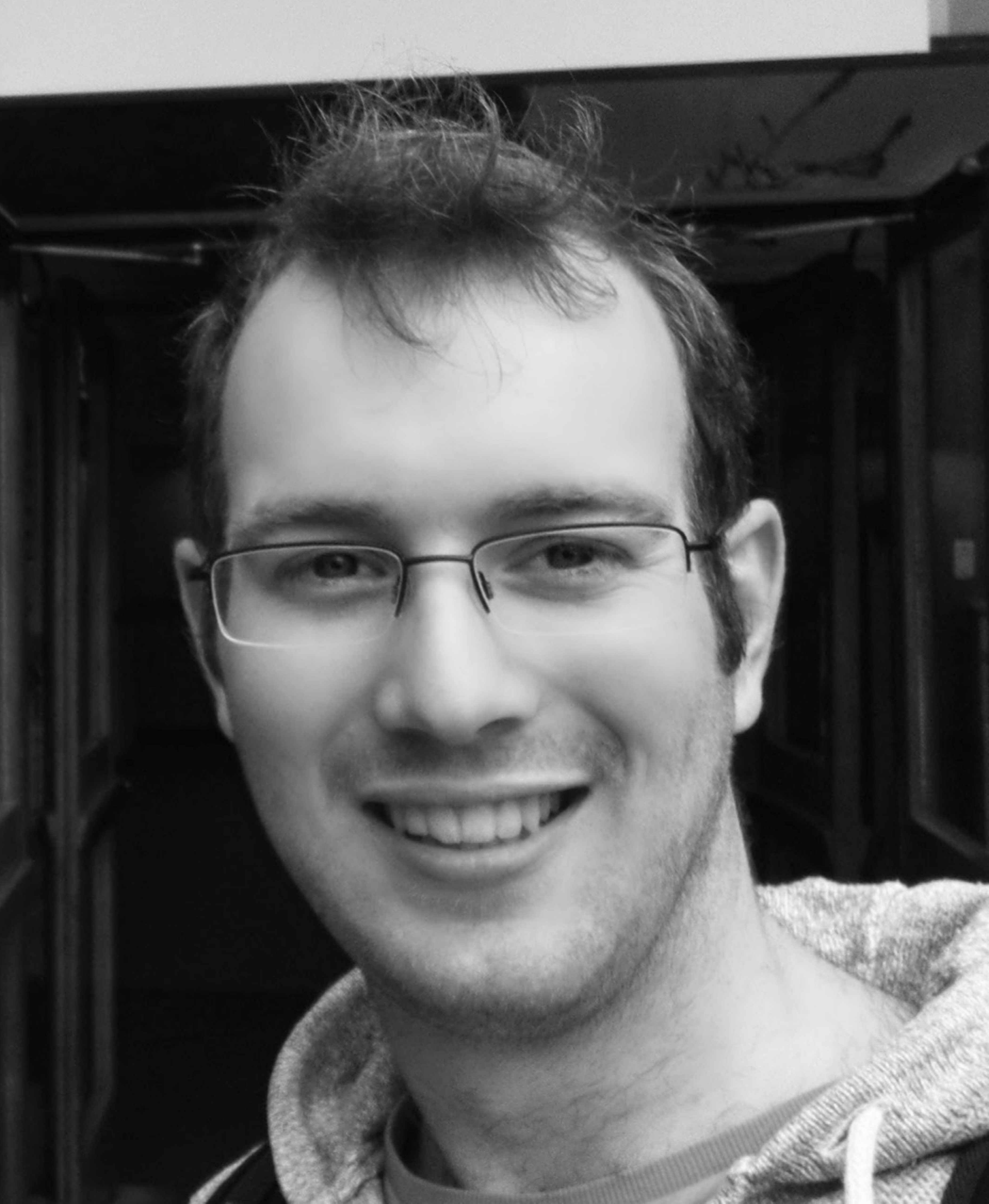}}]{Adam~Herrmann}
	obtained his M.E. degree in Electronic and Computer Engineering in 2017 and his B.Sc. in Engineering Science in 2015, both with first class honours, from University College Dublin. His interests include: Aeronautics, Cyber Security, and Smart Cities, with a focus on Smart Transportation Systems.
\end{IEEEbiography}

\begin{IEEEbiography}[{\includegraphics[width=1in,height=1.25in,clip,keepaspectratio]{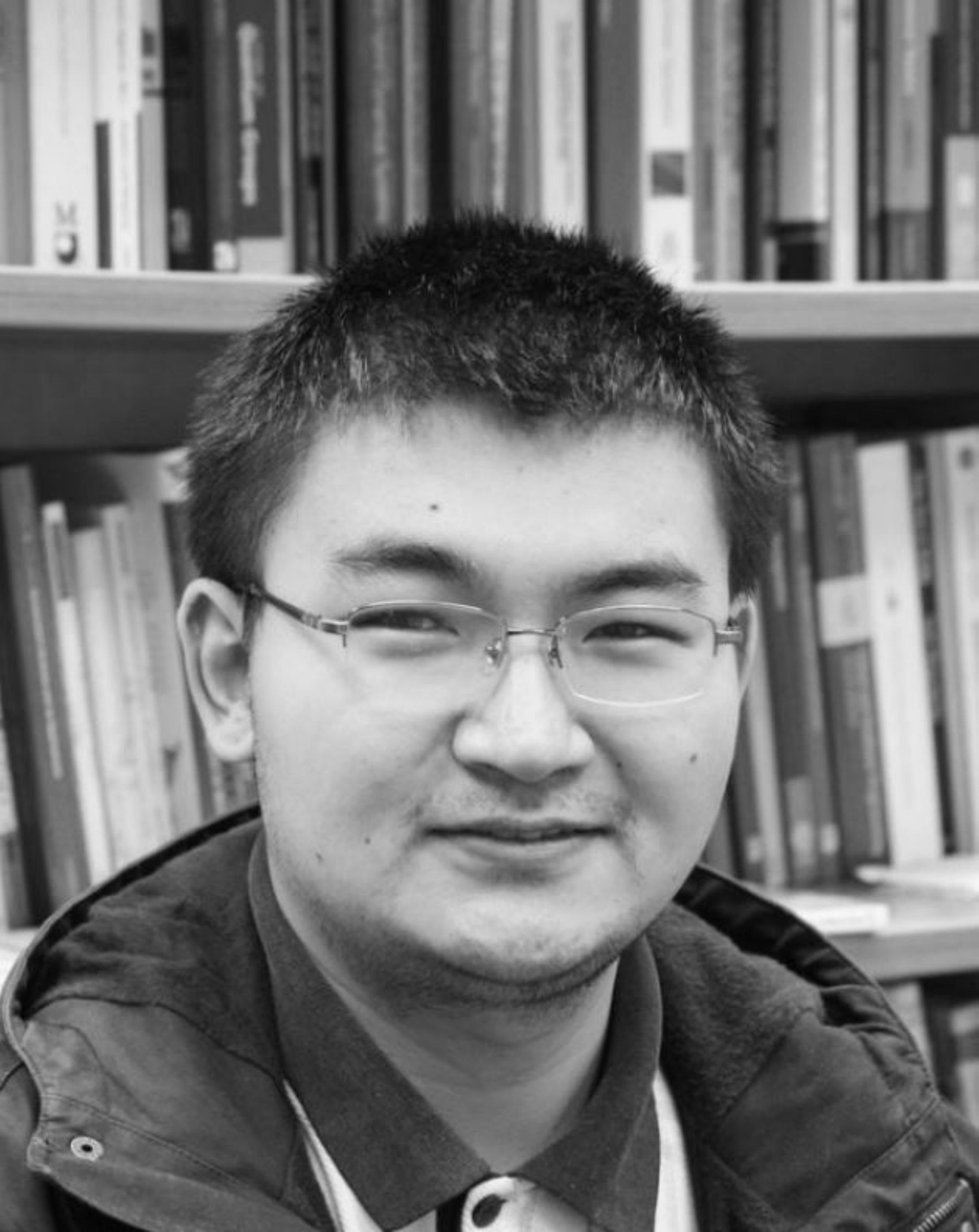}}]{Mingming~Liu}
	received his B.E. degrees in Electronic Engineering with first class honours from National University of Ireland Maynooth in 2011. He obtained his Ph.D. degree from the Hamilton Institute, Maynooth University in 2015. He is currently a post-doctoral research fellow in University College Dublin, working with Prof. Robert Shorten. His current research interests are nonlinear system dynamics, distributed control techniques, modelling and optimisation in the context of smart grid and smart transportation systems. 		
\end{IEEEbiography}

\begin{IEEEbiography}[{\includegraphics[width=1in,height=1.25in,clip,keepaspectratio]{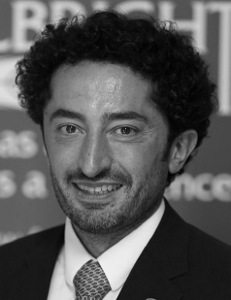}}]{Francesco~Pilla's}
	area of expertise is geospatial analysis and modelling of environmental dynamics, which involve the development of environmental pollution models (air, noise, water) and decision support tools using a GIS platform, in order to facilitate the interoperability of input data and research outcomes between the client/final user and the research team. His work focuses on understanding complex environmental phenomena in order to preempt the impacts resulting from interactions between the human population and the environment.  Before joining UCD, he was at Trinity College Dublin and MIT.\end{IEEEbiography}

\begin{IEEEbiography}[{\includegraphics[width=1in,height=1.25in,clip,keepaspectratio]{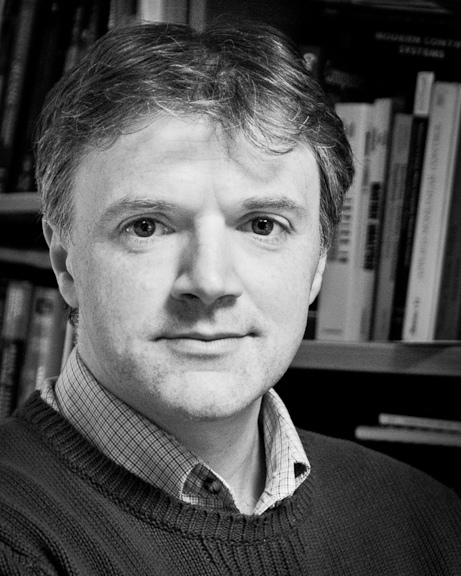}}]{Robert~Shorten}
 is Professor of Control Engineering and Decision Science at UCD,  and retains a part-time appointment with IBM Research. 
 He is a co-founder of the Hamilton Institute at Maynooth University, has been a visiting professor at TU Berlin, and led the {\em Optimisation and 
 Control} team at IBM Research at their Smart Cities Research Lab in Dublin. He has also been a research visitor at Yale University and 
 Technion. His research interests include: Smart Mobility and Smart Cities; Control Theory and Dynamics; Hybrid Dynamical Systems; 
 Networking; Linear Algebra; Collaborative Consumption. He is the Irish member of the European Union Control Association (EUCA) assembly,
  a member of the IEEE Control Systems Society Technical Group on Smart Cities, and a member of the IFAC Technical Committees for 
  Automotive Control, and for Discrete Event and Hybrid Systems. He is a co-author of the recently published books {\em 
 AIMD Dynamics and Distributed Resource Allocation (SIAM 2016)} and {\em Electric and Plug-in Vehicle Networks: Optimisation and Control (CRC Press, Taylor and Francis Group, 2017)}.
\end{IEEEbiography}

\end{document}